\documentclass[draft]{article}

\usepackage{amsmath, amsthm, amsfonts}
\usepackage{mathrsfs}
\usepackage{pstricks, pst-node}
\usepackage{subfigure}
\usepackage{tabularx}
\usepackage{calc}

\newtheoremstyle{break}{\topsep}{\topsep}
     {\itshape}%         Body font
     {}%         Indent amount (empty = no indent, \parindent = para indent)
     {\bfseries}% Thm head font
     {.}%        Punctuation after thm head
     {\newline}%     Space after thm head (\newline = linebreak)
     {}%{\thmname{#1}\thmnumber{ #2}\thmnote{ #3}}%         Thm head spec

\theoremstyle{break}  \newtheorem{theorem}{Theorem}
		      \newtheorem{lemma}[theorem]{Lemma}
		      
\theoremstyle{remark} \newtheorem*{remark}{Remark}

\numberwithin{equation}{section}
\renewcommand{\P}[1] {P\!\left(#1\right)}
\newcommand{\E}[1] {\boldsymbol{E}\!\left(#1\right)}
\newcommand{\Var}[1] {\mathit{Var}\!\left(#1\right)}
\newcommand{\Cov}[2] {\mathit{Cov}\!\left(#1,#2\right)}

\newcommand{\Ind}[1] {1_{\left\{#1\right\}}}
\newcommand{\I}[2] {1_{ #1 #2} }  % 1_uv
\newcommand{\ti}[1] {\mathcal{I}_{\mathbf{#1}}}  % Topological index
\newcommand{\convd} {\xrightarrow{\mathscr{D}}}
\newcommand{\Poiss}[1] {\mathscr{P}_{\!#1}}
\newcommand{\dfk}[2] {\delta_{#1}^{(#2)}}
\newcommand{\pdfk}[2] {\left(\delta_{#1}^{(#2)}\right)}
\newcommand{\Eedges} {\boldsymbol{E}|E|}
\renewcommand{\O}[1] {O\!\left(#1\right)}

\renewcommand{\l} {\left}
\renewcommand{\r} {\right}
\newcommand{\nlbr} {\right.\\\left.}
\newcommand{\lands} {\,\land\,}

\begin{document}

%---------------------------------------------------------------------------------------------------------------------------------
\title{The Covariance of Topological Indices that Depend on the Degree of a Vertex}
\author{Boris Hollas\\[2ex]Theoretische Informatik, Universit\"at Ulm, D-89081 Ulm\\{\small E-mail: hollas@informatik.uni-ulm.de}}
\date{\today}
%---------------------------------------------------------------------------------------------------------------------------------
\maketitle

\begin{abstract}
  We consider topological indices $\mathcal{I}$ that are sums of $f(\deg(u)) f(\deg(v))$, where $\{u,v\}$ are adjacent vertices and $f$ is
  a function. The Randi{\'c} connectivity index or the Zagreb group index are examples for indices of this kind. In earlier work on
  topological indices that are sums of independent random variables, we identified the correlation between $\mathcal{I}$ and the edge set
  of the molecular graph as the main cause for correlated indices. We prove a necessary and sufficient condition for $\mathcal{I}$ having
  zero covariance with the edge set.
\end{abstract}

\renewcommand{\wlog} {without loss of generality}

%=====================
\section{Introduction}
%=====================

For quite some time it has been known that topological indices (graph invariants on molecular graphs) exhibit considerable mutual
correlation \cite{Mot+_TIIRC, Bas+_TINMR}. This is a major problem when performing structure-activity studies as the employed statistical
methods may fail or give little meaningful results on sets of correlated data. Also, strong correlations among a set of topological
indices raise doubt whether these indices describe different and meaningful biological, chemical or physical properties of molecules.

In an attempt to investigate the reasons for these correlations, we used random graphs \cite{Bol_MGT} as a model for chemical graphs and
for topological indices of the form
$$ \ti X (G) = \frac{1}{2} \sum_{\{u,v\} \in E} X_u X_v $$
where $E$ is the edge set of the molecular graph $G=(V,E)$ and $\{ X_v \mid v \in V\}$ is a set of independent random variables with a
common expectation $\E{X}$ \cite{Hollas_CPACD, Hollas_AACD, Hollas_CDBD}. We proved that $\ti X$, $\ti Y$, and $\ti 1$ are \emph{linearly
dependent} for \emph{independent} vertex properties $X, Y$ with $\E X, \E Y >0$ as the number of vertices tends to infinity. For $\E X =
\E Y =0$ however these indices are uncorrelated. Here, $\ti 1$ denotes a topological index with $X_v =1$ for all $v \in V$.

While the random graph model we used in \cite{Hollas_CDBD} encompasses graphs of arbitrary structure, including chemical graphs, the
notion of vertex (or atom) properties $X_v$ that are \emph{independent} of the molecular graph is a serious abstraction from computational
chemistry where atom properties used for topological indices are a \emph{function} of the graph or even the molecule.

In this paper, we use a slightly more general random graph model than the one used in \cite{Hollas_CPACD}. In particular, we consider
graphs on $n$ vertices whose edges are chosen independently with a probability proportional to $1/n$. The latter ensures that the expected
number of edges increases linearly in the number of vertices. We use this to model an approximately linear relation of bonds to vertices
present in molecules. For example, homologous series of aliphatic or aromatic hydrocarbons with $n$ atoms contain $n+c$ bonds for some
constant $c$. Polyphenols contain $\frac{7}{6} n+c$ bonds as each monomer adds 6 atoms and 7 bonds. On the other hand, there is some
variation in the number of bonds for a given number of atoms in a heterogenous set of molecules, which is also true for the random graph
model.

As a more significant difference we consider the vertex properties $X_v$ to be a function of the vertex degree instead of being
independent. Thus, our results are valid for important topological indices such as the Randi\'c connectivity index or Zagreb group
index. We will focus on the crucial covariance between $\ti X$ and $\ti 1$.

%======================
\section{Preliminaries}
%======================

First, we describe the random graph model. For a graph $(V,E)$ let
$$ \I uv = 1_{ \{ \{u,v\} \in E \}} = \begin{cases}1 & \text{if $\{u,v\} \in E$} \\ 0 & \text{else} \end{cases}$$
be the indicator function for $\{ \{u,v\} \in E \}$. For $V=\{1,\dots,n\}$ let $\I uv \; (u,v \in V)$ be
independent random variables with $\P{\I uv = 1} = p$. The space of random graphs $\mathscr{G}(n,p)$ can be identified with the
distribution of $(\I uv)_{u,v \in V}$. We set $p=\alpha/n$ for a fixed parameter $\alpha > 0$ so that $\Eedges =
\binom{n}{2} p \sim \frac{\alpha}{2} n$ as motivated in the introduction.

To describe the vertex properties, let $f: \mathbb{N}_0\to \mathbb{R}$ be a function with $f(0)=0$. We consider the topological index
\begin{equation} \label{eq:tiX}
  \ti X = \ti X (G) = \frac{1}{2} \sum_{ \{u,v\} \in E(G)} X_u X_v
\end{equation}
with
$$ X_v = f(\deg(v)) $$
being the vertex properties and $G \in \mathscr{G}(n,\alpha/n)$ is a random graph. Thus, $f(0)=0$ accounts for isolated vertices being
ignored. Using indicators this can be written as
\begin{equation} \label{eq:tiXI}
  \ti X = \sum_{u<v} X_u X_v \I uv
\end{equation}
which is better suited to employ the expectation operator.

We us the following \emph{notations} throughout the text:

\smallskip
\begin{tabularx}{\textwidth - 1cm}{ccX}
  $\O{f}$           &  \quad denotes \quad\quad &  a function $g$ with $g(x) \leq c f(x)$ for all large $x$ and some constant $c>0$\\
  $X_n \convd X$    &  \quad denotes \quad\quad &  that random variable $X_n$ converges to $X$ in distribution\\
  $a_n \nearrow a$  &  \quad denotes \quad\quad	&  that sequence $(a_n)$ is monotonically increasing and converges to $a$
\end{tabularx}

%====================================
\section{Expectations and Covariance} \label{sec:Exp_Cov}
%====================================

To determine expectation values, we have to eliminate the dependence among $X_u$ and $X_v$ in \eqref{eq:tiXI}. This is achieved by
conditioning for $\{\I uv = 1\}$. If the edge $\{u,v\}$ exists then the degree of $u$ has no effect on the degree of $v$ and vice versa:

\begin{lemma} \label{lem:Iuv}
%------------
  Suppose $u<v$. Then the random variables $(\I u{u'})_{u'>u}$ and $(\I v{v'})_{v'>v}$ are independent with respect to the
  probability measure $\P{\cdot \mid \I uv = 1}$. The same claim holds for $\deg(u)$ and $\deg(v)$.
\end{lemma}
\begin{figure}[h]
  \centering
  \begin{pspicture}(-0.5,-0.5)(3,1)
    \dotnode(1,0.5){U}\nput{d}{U}{u}
    \dotnode(2,0.5){V}\nput{d}{V}{v}
    \dotnode(0,0){U1}\nput{d}{U1}{u'}
    \dotnode(0,1){U2}
    \dotnode(3,0){V1}\nput{d}{V1}{v'}
    \dotnode(3,1){V2}
    \ncline{U}{V}
    \ncline[linestyle=dotted]{U}{U1}
    \ncline[linestyle=dotted]{U}{U2}
    \ncline[linestyle=dotted]{V}{V1}
    \ncline[linestyle=dotted]{V}{V2}
  \end{pspicture}  \caption{We fix edge $\{u,v\}$}
\end{figure}
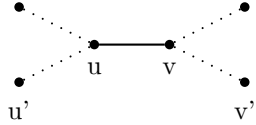
\begin{proof}
  Let $a_{uu'}, a_{vv'} \in \{0,1\}$ for $u'>u$, $v'>v$ and $a_{uv}=1$. We check that for $(a_{uu'})_{u'>u}$, $(a_{uu'})_{u'>u}$ holds
  \begin{multline*}
  \P{ (\I u{u'})_{u'>u} = (a_{uu'})_{u'>u} \lands (\I v{v'})_{v'>v} = (a_{vv'})_{v'>v} \mid \I uv = 1} \\
  = \frac{\P{ (\I u{u'})_{u'>u, u'\neq v} = (a_{uu'})_{u'>u} \lands (\I v{v'})_{v'>v} = (a_{vv'})_{v'>v} \lands \I uv = 1}} {\P{\I uv = 1}} \\
  = \P{ (\I u{u'})_{u'>u, u'\neq v} = (a_{uu'})_{u'>u} \lands (\I v{v'})_{v'>v} = (a_{vv'})_{v'>v}} \\
  = \P{ (\I u{u'})_{u'>u, u'\neq v} = (a_{uu'})_{u'>u} } \P{ (\I v{v'})_{v'>v} = (a_{vv'})_{v'>v}} \\
  = \P{ (\I u{u'})_{u'>u} = (a_{uu'})_{u'>u} \mid \I uv = 1 } \\
                 \cdot \P{ (\I v{v'})_{v'>v} = (a_{vv'})_{v'>v} \mid \I uv = 1}
  \end{multline*}
  If $a_{uv}=0$, both sides are zero. The second claim is a consequence of $\deg(u)$ or $\deg(v)$ being functions of $\I u{u'}$ or $\I
  v{v'}$, respectively.
\end{proof}

We are going to apply lemma \ref{lem:Iuv} to conditional expectations. This motivates the definition
\begin{equation} \label{eq:def_dfk}
  \dfk fk = \E{X_1 \mid \I 12 \I 13 \cdots \I 1{k+1} = 1} , \quad k > 0
\end{equation}
We shall see later why we also need $k>1$. For symmetry reasons, this could as well be defined for a vertex $v \neq 1$ and any set of
distinct vertices $\{u_2,\dots,u_k\}$ different from $v$. As we shall see in section \ref{sec:poissdfk}, $\lim_{n \to \infty} \dfk f k$
exists and is a function of $\alpha$ if $f$ satisfies a condition. Thus, we may regard $\dfk f k$ as almost constant for large $n$.

\begin{lemma} \label{lem:EtiX}
%------------
  $\E{\ti X} = \pdfk f1 ^2 \Eedges$ 
\end{lemma}
\begin{proof}
  \begin{align*}
  \E{\ti X} &= \sum_{u<v} \E{X_u X_v \mid \I uv =1} p                    && \text{by \eqref{eq:tiXI}}\\
            &= \sum_{u<v} \E{X_u \mid \I uv =1} \E{X_v \mid \I uv =1} p  && \text{by lemma \ref{lem:Iuv}}\\
            &= \pdfk f1 ^2 \Eedges                                       && \text{by \eqref{eq:def_dfk}}
  \end{align*}
\end{proof}

\begin{lemma} \label{lem:EtiXti1}
%------------
  $$ \E{\ti X \ti 1} = \l[ \pdfk f1 ^2 \binom{n-2}{2} p + 2 \dfk f1 \dfk f2  (n-2)p + \pdfk f1 ^2 \r] \Eedges $$
\end{lemma}
\begin{proof}
  To dissect the sum
  $$ \E{\ti X \ti 1} = \sum_{u<v} \sum_{u'<v'} \E{X_u X_v \I uv \I {u'}{v'}} $$
  according to $| \{u,v\} \cap \{u',v'\} |$, consider
  $$ S_k = \{ (u,v,u',v') \mid u<v \lands u'<v' \lands |\{u,v\} \cap \{u',v'\}|=k \}, \quad 0 \leq k \leq 2 $$
  Then
  \begin{align}
    |S_0| &= \binom{n}{2} \binom{n-2}{2}  \label{eq:S0}\\
    |S_1| &= 6 \binom{n}{3}               \label{eq:S1}\\
    |S_2| &= \binom{n}{2}                 \label{eq:S2}
  \end{align}
  \eqref{eq:S0} and \eqref{eq:S2} are obvious. To verify \eqref{eq:S1} let $(u,v,u',v') \in S_1$. Exactly two numbers are equal as
  indicated in figure \ref{fig:posS1}. Cases (a), (b) allow just one way to distribute three distinct numbers on $u,v,u',v'$ while there
  are two ways for cases (c), (d).
  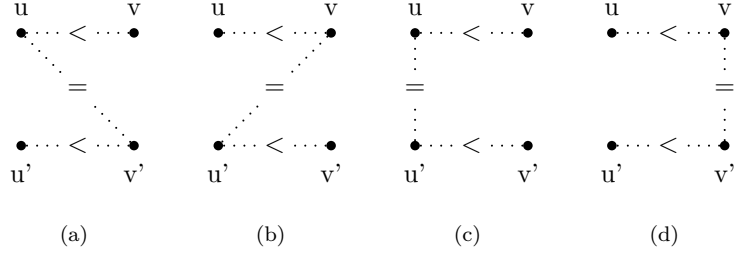
\begin{figure}[h]
    \centering
    \subfigure[]{\begin{pspicture}(-0.5,-0.5)(2,2)
      \dotnode(0,0){U1}\nput{d}{U1}{u'}
      \dotnode(1.5,0){V1}\nput{d}{V1}{v'}
      \dotnode(0,1.5){U}\nput{u}{U}{u}
      \dotnode(1.5,1.5){V}\nput{u}{V}{v}
      \ncline[linestyle=dotted]{U}{V1}\ncput*{$=$}
      \ncline[linestyle=dotted]{U}{V}\ncput*{$<$}
      \ncline[linestyle=dotted]{U1}{V1}\ncput*{$<$}
    \end{pspicture}}
    \subfigure[]{\begin{pspicture}(-0.5,-0.5)(2,2)
      \dotnode(0,0){U1}\nput{d}{U1}{u'}
      \dotnode(1.5,0){V1}\nput{d}{V1}{v'}
      \dotnode(0,1.5){U}\nput{u}{U}{u}
      \dotnode(1.5,1.5){V}\nput{u}{V}{v}
      \ncline[linestyle=dotted]{U1}{V}\ncput*{$=$}
      \ncline[linestyle=dotted]{U}{V}\ncput*{$<$}
      \ncline[linestyle=dotted]{U1}{V1}\ncput*{$<$}
    \end{pspicture}}
    \subfigure[]{\begin{pspicture}(-0.5,-0.5)(2,2)
      \dotnode(0,0){U1}\nput{d}{U1}{u'}
      \dotnode(1.5,0){V1}\nput{d}{V1}{v'}
      \dotnode(0,1.5){U}\nput{u}{U}{u}
      \dotnode(1.5,1.5){V}\nput{u}{V}{v}
      \ncline[linestyle=dotted]{U}{U1}\ncput*{$=$}
      \ncline[linestyle=dotted]{U}{V}\ncput*{$<$}
      \ncline[linestyle=dotted]{U1}{V1}\ncput*{$<$}
    \end{pspicture}}
    \subfigure[]{\begin{pspicture}(-0.5,-0.5)(2,2)
      \dotnode(0,0){U1}\nput{d}{U1}{u'}
      \dotnode(1.5,0){V1}\nput{d}{V1}{v'}
      \dotnode(0,1.5){U}\nput{u}{U}{u}
      \dotnode(1.5,1.5){V}\nput{u}{V}{v}
      \ncline[linestyle=dotted]{V}{V1}\ncput*{$=$}
      \ncline[linestyle=dotted]{U}{V}\ncput*{$<$}
      \ncline[linestyle=dotted]{U1}{V1}\ncput*{$<$}
    \end{pspicture}}
    \caption{Possibilities for $(u,v,u',v') \in S_1$} \label{fig:posS1}
  \end{figure}
  For symmetry reasons, $\E{X_u X_v \I uv \I {u'}{v'}} = \E{X_1 X_2 \I 12 \I 13}$ for all $(u,v,u',v') \in S_1$. Hence, we get
  \begin{equation} \label{eq:EtiXti1}
  \begin{split}
    \E{\ti X \ti 1} &=  |S_0| \E{X_1 X_2 \I 12 \I 34} \\
                    &\quad+  |S_1| \E{X_1 X_2 \I 12 \I 13} \\
		    &\quad+  |S_2| \E{X_1 X_2 \I 12 ^2} \\
                    &=  |S_0| \E{X_1 X_2 \mid \I 12 =1} p^2 \\
                    &\quad+  |S_1| \E{X_1 X_2 \I 13 \mid \I 12 =1} p \\
		    &\quad+  |S_2| \E{X_1 X_2 \mid \I 12 =1} p \\
		    &=  |S_0| \E{X_1 \mid \I 12 =1} \E{X_2 \mid \I 12 =1} p^2 \\
		    &\quad+  |S_1| \E{X_1 \I 13 \mid \I 12 =1} \E{X_2 \mid \I 12 =1} p \\
		    &\quad+  |S_2| \E{X_1 \mid \I 12 =1} \E{X_2 \mid \I 12 =1} p
  \end{split}
  \end{equation}
  by lemma \eqref{lem:Iuv}. With 
  $$ \E{X_1 \I 13 \mid \I 12 =1} = 1/p \E{X_1 \I 12 \I 13} = \dfk f 2 p$$
  $$ \binom{n}{3} = \binom{n}{2} \frac{n-2}{3} $$
  and \eqref{eq:S0}-\eqref{eq:S2}, \eqref{eq:EtiXti1}, we have
  \begin{align*}
    \E{\ti X \ti 1} &=  \pdfk f 1 ^2 \Eedges \binom{n-2}{2} p \\
                    &+  2 \dfk f 1 \dfk f 2 \Eedges (n-2)p \\
		    &+  \pdfk f 1 ^2 \Eedges
  \end{align*}
\end{proof}

\begin{remark}
  With $f \equiv 1$, lemma \ref{lem:EtiX} and the help of Mathematica follows $\Var{\ti 1} = \Eedges (1-p)$, as it should be.
\end{remark}

We combine the results of this section in

\begin{theorem} \label{thm:Cov}
%--------------
  If $\dfk f 1$, $\dfk f 2$ are bounded in $n$ then
  $$ \Cov{\ti X}{\ti 1} = \begin{cases}
    0   &   \text{if $\dfk f 1 = 0$} \\
    \l[ \pdfk f 1 ^2 \l( 1+2\alpha \l(\frac{\dfk f 2}{\dfk f 1} - 1\r) \r) + \O{\frac{1}{n}} \r] \Eedges   &   \text{else}
  \end{cases} $$
\end{theorem}
\begin{proof}
  By lemma \ref{lem:EtiX} and lemma \ref{lem:EtiXti1},
  \begin{multline*}
    \Cov{\ti X}{\ti 1} = \l[ \pdfk f 1 ^2 \binom{n-2}{2} p +  2 \dfk f 1 \dfk f 2 (n-2)p \nlbr
		     + \pdfk f 1 ^2  - \pdfk f 1 ^2 \Eedges \r] \Eedges
  \end{multline*}
  Using $\binom{n-2}{2} - \binom{n}{2} = 3-2n$, this can be written as
  \begin{multline*}
    \Cov{\ti X}{\ti 1} = \l[\pdfk f 1 ^2 (1+(3-2n)p) +  2 \dfk f 1 \dfk f 2 (n-2)p \r] \Eedges \\
    = \begin{cases}
        0   &   \text{if $\dfk f 1 = 0$} \\
	\l[\pdfk f 1 ^2 \l( 1 + 2p\l(\frac{\dfk f 2}{\dfk f 1} n - n\r) \r) +  \pdfk f 1 ^2 p \l(3 - 2\frac{\dfk f 2}{\dfk f 1}\r) \r] &\\
												\hfill\cdot\Eedges   &   \text{else}
   \end{cases}
  \end{multline*}
  The assertion follows with $p = \alpha/n$.
\end{proof}

\begin{remark}
  We will prove in theorem \ref{thm:dfk} in section \ref{sec:poissdfk} that all $\dfk f k$ are in fact bounded in $n$ if $f \in \O{x}$.
\end{remark}

Yet, it is not clear  whether $\Cov{\ti X}{\ti 1} \neq 0$ for $\dfk f 1 \neq 0$. This is dealt with in the next section.

%=================================================
\section{The Poisson Distribution and  $\dfk f k$} \label{sec:poissdfk}
%=================================================

For the proof of the following theorem recall that for random variables $X_n, X$ holds
$$ X_n \convd X$$
iff
$$ \E{f(X_n)} \to \E{f(X)} $$
for all bounded and continuous functions $f: \mathbb{R} \to \mathbb{R}$. This does not hold for arbitrary unbounded functions
$f$. Therefore, we require that $f \in \O{x}$ in this section. While this does not seem to be the most general restriction, it facilitates
the following elaborations.

\begin{theorem} \label{thm:dfk}
%--------------
  For all $f \in \O{x}$ and all $k$ holds
  $$ \lim_{n \to \infty} \dfk f k = \E{f(k+\Poiss{\alpha})} = \sum_{j=0}^\infty f(k+j) \frac{\alpha^j}{j!} e^{-\alpha} $$ where
  $\Poiss{\alpha}$ is the Poisson distribution with parameter $\alpha$.
\end{theorem}
\begin{proof}
  By definition \eqref{eq:def_dfk},
  \begin{align}
    \dfk f k &= \E{f(\deg(1)) \mid \I 12 \I 13 \cdots \I 1{k+1} = 1} \notag\\
	     &= \E{f \bigg( k+\sum_{j=k+2}^n \I 1j \bigg) }      \label{eq:dfkE}
  \end{align}
  Since $p=\alpha/n$, Poisson's limit theorem gives
  $$ \sum_{j=k+2}^n \I 1j \convd \Poiss{\alpha} \quad (n \to \infty) $$ 
  The function $f: \mathbb{N}_0 \to \mathbb{R}$ can be extended to a continuous function $f: \mathbb{R} \to \mathbb{R}$ in an arbitrary
  way. Hence, the continuity theorem gives
  $$ f \l( k+\sum_{j=k+2}^n \I 1j \r) \convd f(k+\Poiss{\alpha}) \quad (n \to \infty) $$
  For all bounded and continuous functions $f^*$ follows by \eqref{eq:dfkE}
  \begin{equation}  \label{eq:df*kE}
     \dfk {f^*} k \to \E{f^*(k+\Poiss{\alpha})} \quad (n \to \infty)
  \end{equation}
  If $f$ is also bounded the claim follows. To prove \eqref{eq:df*kE} for unbounded $f$ we cut $f$ off above a limit to divide $f$ into a
  bounded and an unbounded part. We show that the unbounded part tends to zero as the limit tends to infinity.

  To begin with, let $|f(x)|\leq x$ for all $x$ and let $f$ be unbounded. Then there is a sequence of integers $(m_l)$ such that \wlog
  $f(m_l) \nearrow \infty$ for $l \to \infty$ and $f(m_l) >0$ for all $l$. Let be
  $$ c_m(x) = \begin{cases} x & \text{if } |x|<m \\ 
			    0 & \text{else}
	      \end{cases} $$
  and
  $$ \tilde{c}_m(x) = \begin{cases} 0 & \text{if } |x|<m \\ 
			    x & \text{else}
	      \end{cases} $$
  Let $S_n := k+\sum_{j=k+2}^n \I 1j$. Then
  \begin{align*}
    |\E{ (\tilde{c}_{m_l} \circ f)(S_n) }| &=    \l|\E{ f(S_n) \Ind{f(S_n) > m_l} }\r|      \\
				  &\leq \E{ S_n \Ind{f(S_n) > f(m_l)} }
  \intertext{since $0 \leq f(m_l) \leq m_l$}
				  &= \E{ S_n \Ind{S_n > m_l} }	     
  \intertext{since $f(m_l)$ increases monotonically}
				  &<	\E{ S_n \frac{S_n}{m_l} }		     \\
				  &=	\frac{1}{m_l} \l[ \Var{S_n}+(\E{S_n})^2 \r]     \\
				  &=	\frac{1}{m_l} \l[ \O{n}p(1-p) + (\O{n}p)^2 \r]  \\
				  &=  \O{1/m_l}
  \end{align*}
  By linearity of expectation follows
  \begin{equation}  \label{eq:EcmfS_n}
    \E{ (\tilde{c}_{m_l} \circ f)(S_n) }  = \O{1/m_l}
  \end{equation}
  for all $f\in\O{x}$. Thus,
  \begin{align*}
    \lim_{n\to\infty} \dfk f k  &= \lim_{l\to\infty} \lim_{n\to\infty} \E{f(S_n)} \\
  \intertext{by \eqref{eq:dfkE}}
			         &= \lim_{l\to\infty} \lim_{n\to\infty} [ \E{(c_{m_l} \circ f)(S_n)} + \E{(\tilde{c}_{m_l} \circ f)(S_n)} ]\\ 
				 &= \lim_{l\to\infty} [ \E{(c_{m_l} \circ f)(k+\Poiss{\alpha})} + \O{1/m_l} ]
  \intertext{by \eqref{eq:df*kE} and \eqref{eq:EcmfS_n}}
				 &= \E{f(k+\Poiss{\alpha})}
  \end{align*}
  by the convergence theorem of Lebesgue.
\end{proof}

With the help of theorem \ref{thm:dfk} we are able to answer the question raised at the end of section \ref{sec:Exp_Cov}:

\begin{theorem} \label{thm:Cov0}
  For $n\to\infty$ and $f \in \O{x}$ holds: $\ti X$ and $\ti 1$ have covariance zero if and only if $\lim_{n\to\infty} \dfk f 1 =0$.
\end{theorem}
%------------
\begin{proof}
  Assume that $\lim_{n\to\infty} \dfk f 1 \neq 0$ and $\lim_{n\to\infty} \Cov{\ti X}{\ti 1}=0$. By theorem \ref{thm:Cov} follows
  $$ \lim_{n\to\infty} \frac{\dfk f 2}{\dfk f 1} = 1 - \frac{2}{\alpha} $$
  With theorem \ref{thm:dfk} we get
  \begin{align*}
    \sum_{j=0}^{\infty} f(2+j) \frac{\alpha^{j}}{j!}  &= \l(1 - \frac{1}{2\alpha}\r) \sum_{j=0}^{\infty} f(1+j) \frac{\alpha^{j}}{j!} \\
						      &= \sum_{j=0}^{\infty} f(1+j) \frac{\alpha^{j}}{j!} - \frac{1}{2} \sum_{j=0}^{\infty} f(1+j) \frac{\alpha^{j-1}}{j!}
  \end{align*}
  We multiply by $\alpha$ and substitute $j$ with $j-1$ in the first two series to get
  $$ \sum_{j=1}^{\infty} f(1+j) \frac{\alpha^{j}}{(j-1)!}  =   \sum_{j=1}^{\infty} f(j) \frac{\alpha^{j}}{(j-1)!} - \frac{1}{2} 
  \sum_{j=0}^{\infty} f(1+j) \frac{\alpha^{j}}{j!} $$
  Hence,
  $$ \frac{1}{2} f(1) \alpha^0 + \sum_{j=1}^{\infty} \l[ f(1+j) \l(1+\frac{1}{2j}\r) - f(j) \r] \frac{\alpha^j}{(j-1)!} = 0 $$ By theorem
  \ref{thm:dfk}, this series converges for all $\alpha >0$. By the identity theorem for power series follows that all coefficients are
  zero. By induction thus follows $f \equiv 0$, which contradicts $\lim_{n\to\infty} \dfk f 1 \neq 0$.

  The opposite direction follows by theorem \ref{thm:Cov} and theorem \ref{thm:dfk}.
\end{proof}

%===================
\section{Discussion}
%===================

We have seen that $\dfk f 1$ is an important quantity for the covariance of the topological indices we consider. Theorem \ref{thm:dfk}
shows that $\dfk f k$ does not depend on $n$ for large $n$. This justifies definition \eqref{eq:def_dfk} since we do not want $\dfk f k$
to be very different for graphs of different size. Also, theorem \ref{thm:dfk} provides a way to approximately compute $\dfk f k$. If we
substitute $X_v$ by $X_v - \dfk f 1$ in \eqref{eq:tiX}, the resulting index is uncorrelated to $\ti 1$.

As a drawback, we require $f \in\O{x}$ in section \ref{sec:poissdfk}. Theorem \ref{thm:dfk} may not be valid if $f$ increases very
steeply. However, it should be possible to derive an upper limit similar to \eqref{eq:EcmfS_n} for functions $f$ with a higher rate of
growth than $\O{x}$.

In \cite{Hollas_CDBD}, we proved that topological indices (with independent vertex properties) are necessarily correlated if the vertex
properties have expectations not equal to zero. Theorem \ref{thm:Cov0} does not give this result as it is an assertion on covariance
only. The next step will therefore be an examination of correlations within this setting.

\bibliographystyle{unsrt}
\bibliography{/home/boris/papers/bib/descr,/home/boris/papers/bib/hollas,/home/boris/papers/bib/math}

\end{document}